\documentclass[11pt]{article}
\usepackage{amsfonts,amsmath,amssymb}
\usepackage{hyperref}
\usepackage{graphicx}
\usepackage{color}

\newcommand{\halmos}{\rule{1ex}{1.4ex}}
\def \qed {\nopagebreak{\hspace*{\fill}$\halmos$\medskip}}

\newtheorem{theorem}{Theorem}[section]
\newtheorem{proposition}[theorem]{Proposition}
\newtheorem{corollary}[theorem]{Corollary}
\newtheorem{conjecture}[theorem]{Conjecture}
\newtheorem{lemma}[theorem]{Lemma}

\newtheorem{remark}[theorem]{Remark}
\newtheorem{defi}[theorem]{Definition}

\newcommand{\bt}{\begin{theorem}}
\newcommand{\et}{\end{theorem}}
\newcommand{\bl}{\begin{lemma}}
\newcommand{\el}{\end{lemma}}
\newcommand{\bp}{\begin{proposition}}
\newcommand{\ep}{\end{proposition}}
\newcommand{\bcor}{\begin{corollary}}
\newcommand{\ecor}{\end{corollary}}
\newcommand{\br}{\begin{remark}\rm}
\newcommand{\er}{\end{remark}}
\newcommand{\bcon}{\begin{conjecture}}
\newcommand{\econ}{\end{conjecture}}
\newcommand{\bd}{\begin{defi}}
\newcommand{\ed}{\end{defi}}
\newcommand{\bee}{\begin{enumerate}}
\newcommand{\eee}{\end{enumerate}}

\newcommand{\be}{\begin{equation}}
\newcommand{\ee}{\end{equation}}

\newcommand{\Mi}{{\cal M}}

\newcommand{\eps}{\epsilon}

\newcommand{\R}{{\mathbb R}}
\newcommand{\N}{{\mathbb N}}
\newcommand{\Z}{{\mathbb Z}}

\renewcommand{\P}{{\mathbb P}}
\newcommand{\E}{{\mathbb E}}

\newcommand{\asto}[1]{\underset{{#1}\to\infty}{\longrightarrow}}
\newcommand{\Asto}[1]{\underset{{#1}\to\infty}{\Longrightarrow}}

\begin{document}

\makeatletter\@addtoreset{equation}{section}
\makeatother\def\theequation{\thesection.\arabic{equation}}

\renewcommand{\labelenumi}{{(\roman{enumi})}}

\title{One-dimensional Voter Model Interface Revisited}
\author{Siva R.\ Athreya$^{\,1}$ \and Rongfeng Sun$^{\,2}$}

\maketitle

\footnotetext[1]{Indian Statistical Institute, 8th Mile Mysore Road, Bangalore,
560059 India. Email: athreya@isibang.ac.in}

\footnotetext[2]{Dept.\ of Mathematics, National University of Singapore, 10 Lower Kent Ridge Road, 119076 Singapore. Email: matsr@nus.edu.sg}

\begin{abstract}
We consider the voter model on $\Z$, starting with all $1$'s to the left of the origin and all $0$'s to the right of the origin.
It is known that if the associated random walk kernel $p(\cdot)$ has zero mean and a finite $\gamma$-th moment for any $\gamma>3$, then the evolution of the
boundaries of the interface region between 1's and 0's converge in distribution to a standard Brownian motion $(B_t)_{t\geq 0}$ under diffusive scaling of space and time. This convergence fails when $p(\cdot)$ has
an infinite $\gamma$-th moment for any $\gamma<3$, due to the loss of tightness caused by a few isolated 1's appearing deep within the regions of all $0$'s
(and vice versa) at exceptional times. In this note, we show that as long as $p(\cdot)$ has a finite second moment, the measure-valued process
induced by the rescaled voter model configuration is tight, and converges weakly to the measure-valued process $(1_{x<B_t} {\rm d}x)_{t\geq 0}$.

\end{abstract}

\vspace{12pt}

\section{Introduction}
The voter model on $\Z$ is an interacting particle system with state space $\Omega:=\{0,1\}^\Z$. At each time $t\geq 0$, we denote the state of the voter model by
$\eta_t:=(\eta_t(x))_{x\in\Z}\in \Omega$, where $\eta_t(x) \in \{0,1\}$ encodes the opinion of the voter at site $x$ at time $t$. Independently for each pair $x,y\in\Z$,
the opinion at $x$ is replaced by the opinion at $y$ (also called {\em resampling}) with exponential rate $p(y-x)$, where $p(\cdot):=(p(x))_{x\in\Z}$ is the probability distribution of the increments of an irreducible random walk on $\Z$. Formally, the voter model has generator
\be
(Lf)(\eta) = \sum_{x,y\in\Z} p(y-x) (f(\eta^{x,y})-f(\eta)),
\ee
where $\eta\in \Omega$, $\eta^{x,y}(z)=\eta(z)$ for all $z\neq x$ and $\eta^{x,y}(x)=\eta(y)$, and $f: \{0,1\}^\Z\to\R$ depends only on a finite number of coordinates.
For classic results on the voter model, see~\cite{L85}.

We consider the voter model with the heavy-side initial configuration
\be\label{heavyside}
\eta_0(x) = \left\{
\begin{aligned}
1 \qquad \mbox{ if } x \leq 0, \\
0 \qquad \mbox{ if } x \geq 1.
\end{aligned}
\right.
\ee
For each $t\geq 0$, we denote the positions of the leftmost 0 and the rightmost 1 respectively by
$$
l_t := \inf\{x \in \Z: \eta_t(x) = 0\}  \qquad \mbox{and} \qquad r_t := \sup\{x \in \Z: \eta_t(x) =1 \}.
$$
The region between $l_t$ and $r_t$ is called the {\em interface region}, where the voter model configuration $\eta_t$ has a mixture of 0's and 1's. The voter model configuration viewed from the leftmost $0$, i.e., $\tilde \eta_t(x):= \eta_t(l_t+x)$ for $x\in\N$, is called the {\em interface process}, which is a Markov process with state space $\{\tilde\eta\in \{0,1\}^\N : \sum_{x\in\N} \tilde\eta(x)<\infty\}$.

In~\cite{CD95}, Cox and Durrett studied the interface process $\tilde\eta_t$. They observed that $\tilde\eta_t$ is positive-recurrent if and only if the
distribution of the interface size $r_t-l_t$ is tight over times $t\geq 0$, which they verified under the assumption that $p(\cdot)$
has a finite third moment, i.e., $\sum_{x\in\Z} |x|^3p(x)<\infty$. This tightness result was later extended by Belhaouari, Mountford and Valle~\cite{BMV07} to $p(\cdot)$ with a finite second moment, which they showed to be optimal in the sense that tightness is lost if $\sum_{x\in\Z} |x|^\gamma p(x)=\infty$ for some $\gamma<2$. An alternative proof of the tightness of $\{r_t-l_t\}_{t\geq 0}$, under the finite second moment assumption on $p(\cdot)$, was given recently by Sturm and Swart in~\cite{SS08}.

Assume without loss of generality that $\sum_{x\in\Z} xp(x)=0$, and $\sigma^2:=\sum_{x\in\Z} x^2p(x)<\infty$. Cox and Durrett~\cite{CD95} also observed that when $\{r_t-l_t\}_{t\geq 0}$ is tight, the finite-dimensional distributions of
$$
\Big(\frac{l_{tN^2}}{\sigma N}\Big)_{t\geq 0} \qquad \mbox{and} \qquad \Big(\frac{r_{tN^2}}{\sigma N}\Big)_{t\geq 0}
$$
converge to those of a standard Brownian motion $(B_t)_{t\geq 0}$ as $N\to \infty$. It is then natural to ask whether $(l_{tN^2}/\sigma N, r_{tN^2}/\sigma N)_{t\geq 0}$
converges in distribution to $(B_t, B_t)_{t\geq 0}$ in path space, i.e., the product space $D([0,\infty),\R)^2$ where $D([0,\infty), \R)$ is the
space of c\`adl\`ag paths equipped with the Skorohod topology. Such a path level convergence would imply that in the diffusive scaling limit, the interface region becomes sharp uniformly on finite time intervals, and the motion of the interface location converges weakly to a Brownian motion. This was established by Newman, Ravishankar and Sun in~\cite{NRS05} under the assumption that $p(\cdot)$ has a finite fifth moment. It was later extended by Belhaouar et al.~\cite{BMSV06} to all $p(\cdot)$ with a finite $\gamma$-th moment for some $\gamma>3$.

It was also pointed out in~\cite{BMSV06} that if $\sum_{x\in\Z} |x|^\gamma p(x)=\infty$ for some $\gamma<3$, then $(l_{tN^2}/\sigma N, r_{tN^2}/\sigma N)_{t\geq 0}$ loses tightness in path space as $N\to\infty$, because there exist exceptional times when 1's appear deep in the region of all 0's (and vice versa) due to the heavy tail of $p(\cdot)$. However, we expect such 1's (and 0's) to be rare and sparse when they do appear, because $\{r_t-l_t\}_{t\geq 0}$ remains tight as long as $p(\cdot)$ has a finite second moment. If we can suitably discount such rare 1's (and 0's), then we should be able to recover the tightness of $(l_{tN^2}/\sigma N, r_{tN^2}/\sigma N)_{t\geq 0}$ as $N\to\infty$, and hence assert the weak convergence of the interface evolution to a Brownian motion. One way to discount such rare 1's and 0's and to restore path level tightness is by suppressing the resampling of voter model opinions involving sites $x,y\in\Z$ with $|y-x|\geq N^\eps$, for some $\eps>0$ depending on $p(\cdot)$. This was the approach taken in~\cite[Theorem 1.3]{BMSV06}, which requires $p(\cdot)$ to have a finite $\gamma$-th moment for some $\gamma>2$.

In this note, we take an alternative approach to the convergence of the voter model interface evolution, which naturally discounts isolated 1's and 0's, and where finite second moment is the natural assumption on $p(\cdot)$. More precisely, we consider the measure-valued process $(\mu_t)_{t\geq 0}$ induced by the voter model configurations $(\eta_t)_{t\geq 0}$, defined by
\be\label{mut}
\mu_t(\cdot) := \sum_{x\in\Z} \eta_t(x) \delta_x(\cdot).
\ee
The state space of $(\mu_t)_{t\geq 0}$ is $\Mi(\R)$, the space of non-negative Radon measures on $\R$ equipped with the vague topology, so that $\mu_n\to \mu$ in $\Mi(\R)$ if and only if $\int f{\rm d}\mu_n\to\int f{\rm d}\mu$ for all $f\in C_c(\R)$, where $C_c(\R)$ denotes the space of continuous functions with compact support on $\R$. For each $N>1$,
we define the rescaled measure-valued process $\mu^{N}_t$ by
\be\label{muNt}
\int f(x) \mu^{N}_t({\rm d}x) := \frac{1}{N} \int f\Big(\frac{x}{N}\Big) \mu_{tN^2}({\rm d}x) \qquad \mbox{for all } f\in C_c(\R).
\ee
Let $D([0,\infty), \Mi(\R))$ denote the space of right-continuous paths in $\Mi(\R)$ with left-hand limits, equipped with the Skorohod topology.
\medskip

Here is our main result.
\bt\label{T:vmi} Assume that $\sum_{x\in\Z} xp(x)=0$ and $\sigma^2:=\sum_{x\in\Z} x^2p(x)<\infty$. Then the distribution of $(\mu^{N}_t)_{t\geq 0}$ on $D([0,\infty), \Mi(\R))$ converges weakly to that of $(\nu_t)_{t\geq 0}:= (1_{\{x<\sigma B_t\}} {\rm d}x)_{t\geq 0}$ as $N\to\infty$, where $(B_t)_{t\geq 0}$ is a standard Brownian motion.
\et
Theorem~\ref{T:vmi} shows that as long as $p(\cdot)$ has a finite second moment, the voter model interface evolution is tight in the measure-valued sense, and converges weakly to a sharp interface following a Brownian path.
\medskip

One may ask what type of measure-valued processes arise in the scaling limit if we take a sequence of voter model initial configurations $\eta^N_0\in \{0,1\}^\Z$, such that $\mu^N_0$ converges vaguely to a limiting measure $\nu_0({\rm d}x)=f_0(x) {\rm d}x$ for some $f_0: \R\to [0,1]$. The answer is that the limit should be the so-called {\em continuum-sites stepping-stone model with Brownian migration} (CSMBM). See~\cite{Z03, Z08} for the CSMBM on the real line and on the one-dimensional torus, and see the references therein for results on continuum-sites stepping-stone models in general. In~\cite{Z03, Z08}, the distribution of the CSMBM was specified using a finite collection of dual coalescing Brownian motions running backwards in time. With the aid of the so-called {\em Brownian web} (see e.g.~\cite{FINR04}), which constructs simultaneously coalescing Brownian motions starting from every point in space and time, one can in fact give a graphical construction of the CSMBM in the same spirit as the graphical construction of the voter model from the dual family of coalescing random walks (see e.g.~\cite{L85}).
Almost surely, for any time $t>0$, the coalescing Brownian motions in the Brownian web starting from every point in $\R$ at time $t$, running backwards in time, determine an ergodic locally finite point configuration $\cdots <x_i<x_{i+1}<\cdots$ on $\R$, such that all coalescing Brownian motions starting from $(x_i, x_{i+1})$ at time $t$ coalesce into a single point $y_i$ at time $0$, and $y_i<y_{i+1}$ for all $i\in\Z$. The configuration of the CSMBM at time $t$ is then given by $\nu_t({\rm d}x)=f_t(x) {\rm d}x$, where independently for each $i\in\Z$, $f_t=1$ on $(x_i, x_{i+1})$ with probability $f_0(y_i)$, and $f_t=0$ on $(x_i, x_{i+1})$ with probability $1-f_0(y_i)$. Our proof of the tightness of $\{\mu^N_\cdot\}_{N>1}$ in Theorem~\ref{T:vmi} is in fact independent of the initial configuration $\eta_0$, and hence applies in this more general setting as well. Proving convergence of the finite-dimensional distributions of $\{\mu^N_\cdot\}_{N>1}$ however requires more care. We will not work out the details here and instead leave it open for the reader, since our main interest is the convergence of the voter model interface in Theorem~\ref{T:vmi}.
\medskip

We remark that scaling limits of voter models in the measure valued setting have been considered before. In dimension $1$, if the voter model has a long range kernel
with scale $M_N$, space and time are rescaled respectively by $M_N\sqrt{N}$ and $N$, then when $M_N/\sqrt{N}\to \rho$ for some $\rho>0$, the density of 1's is shown in~\cite{MT95} to converge to the solution of an SPDE; when $M_N/\sqrt{N}\to\infty$, the voter model is shown in~\cite[Theorem 1.1]{CDP00} to converge to super-Brownian motion. As mentioned in~\cite{CDP00}, when $M_N/\sqrt{N}\to 0$, the scaling limit of the voter model is believed to be the CSMBM described above. In dimensions $2$ and higher, the voter model has been shown to converge to super-Brownian motion, see e.g.~\cite{CDP00, BCG01}.

\section{Proof}
To prove Theorem \ref{T:vmi}, it suffices to show: (1) tightness of $\{(\mu^{N}_t)_{t\geq 0}\}_{N>1}$ on $D([0,\infty), \Mi(\R))$, which is where the main technical difficulty lies; (2) convergence of the finite dimensional distributions of $(\mu^{N}_t)_{t\geq 0}$ to that of $(\nu_t)_{t\geq 0}:= (1_{\{x<\sigma B_t\}} {\rm d}x)_{t\geq 0}$. We note that our tightness proof is independent of the initial configuration $\eta_0$ of the voter model.

\subsection{Tightness}
By Jakubowski's tightness criterion (see e.g.\ \cite[Theorem 3.6.3]{DA93}),
 $\{(\mu^{N}_t)_{t\geq 0}\}_{N>1}$ is tight on $D([0,\infty), \Mi(\R))$ if the following conditions are satisfied:
\bee
 \item[{\bf (J1)}] (Compact Containment) For each $T >0$ and $\epsilon >0$, there exists a compact set  $K_{T,\epsilon} \subset  \Mi(\R)$ such that for all $N>1$,
\be \label{cc}
\P\big(\mu^{N}_t \in K_{T,\epsilon}, \,\forall \, 0 \leq t \leq T \big)  \geq 1 - \epsilon;
\ee

 \item[{\bf (J2)}] (Tightness of Evaluations) For each $f\in C_c^2(\R)$, the space of twice continuously differentiable real-valued functions on $\R$ with compact support,
 define
 \be\label{XNt}
  X_{t}^{N} := X_{t}^{N} (f) := \int f(x) \mu^{N}_t ({\rm d}x) = \frac{1}{N} \sum_{x \in \Z} f\Big(\frac{x}{N}\Big) \eta_{tN^{2}}(x).
  \ee
  Then $\{(X^N_t)_{t\geq 0}\}_{N>1}$ is tight on $D([0,\infty), \R)$.
\eee
Condition {\bf (J1)} is easily seen to hold, because for each $N>1$ and $t\geq 0$,
$$
\mu^N_t([-m,m]) =   \frac{1}{N}\sum_{x\in \Z \cap [-mN, mN]} \eta_{tN^{2}}(x) \leq 2m+1 \quad \mbox{ for all } m\in\N,
$$
and $K:=\{\nu \in \Mi(\R) : \nu([-m, m]) \leq 2m+1  \ \forall\, m \in \N\}$ is a compact subset of $\Mi(\R)$.

We will verify {\bf (J2)} by verifying Aldous' tightness criterion (see e.g.\ \cite[Theorem 3.6.4]{DA93}) for $\{(X^N_t)_{t\geq 0}\}_{N>1}$ in
$D([0,\infty), \R)$, which reduces to the following conditions:
\bee
\item[{\bf (A1)}] For each rational  $t \geq 0$,  $\{X^{N}_{t}\}_{N>1}$ is tight in $\R$;
\item[{\bf (A2)}] For $T>0$, let $\tau_{N}$ be stopping times bounded by $T$, and let $\delta_{N} \downarrow 0$ as $N \rightarrow \infty$.
Then
$$ \lim_{N \to \infty} \P( | X^{N}_{\tau_{N} + \delta_{N}} - X^{N}_{\tau_{N}} | > \epsilon) = 0.$$
\eee

\begin{remark}{\rm
Before embarking on the verification of {\bf (A1)--(A2)}, we first briefly explain the main technical difficulty in proving tightness. Typically, Aldous' criterion is verified by verifying a criterion of Joffe and M\'etivier (see e.g.\ \cite[Theorem 3.6.6]{DA93}), which requires bounds on the so-called local coefficients
of first and second order, given respectively by
$$
\alpha^N_t := L^N X^N_t  \qquad \mbox{and} \qquad \beta^N_t := L^N((X^N_t)^2)-2X^N_t \alpha^N_t,
$$
where $L^N$ is the generator for the measure-valued process $\mu^N_t$. The coefficients $\alpha^N_t$ and $\beta^N_t$ encode the drift and quadratic variation
of $(X^N_t)_{t\geq 0}$. For purposes of illustration, let us consider $\mu^N_t$ obtained by diffusively rescaling the voter model $\eta^N_t$ on the torus $\Z/(2N\Z)$, identified with $[-N+1, N]\cap\Z$, with $\eta^N_0(\cdot)= 0$ on $[-N+1,0]$ and $\eta^N_0(\cdot)= 0$ on $[1,N]$. Assume that $p(\cdot)$ is symmetric, and its projection $p^N(\cdot)$ on the torus is used to define $(\eta^N_t)_{t\geq 0}$. Take $f\equiv 1$ on the continuum torus $[-1,1]$ with $-1$ identified with $1$. Then
$$
X^N_t=X^N_t(f)=\frac{1}{N}\sum_{x\in [-N+1,N]\cap\Z} \eta^N_{tN^2}(x)
$$
is a martingale, and hence $\alpha^N_t= 0$ for all $t\geq 0$. A simple calculation shows that
$$
\beta^N_t = \frac{1}{2} \sum_{x,y\in [-N+1,N]\cap\Z} p^N(y-x) 1_{\{\eta^N_{tN^2}(x)\neq \eta^N_{tN^2}(y)\}},
$$
which will be $O(1)$ only if $\eta^N_{tN^2}$ segregates into $O(1)$ number of intervals with mostly all 1's or all 0's on each interval. Establishing such segregation
of 0's and 1's is the main difficulty in proving tightness, which Joffe and Metivier's criterion does not help to simplify. Instead, we will proceed by
a direct verification of Aldous' criterion, using the duality between the voter model and coalescing random walks.
}\end{remark}
\bigskip

\noindent
{\bf Proof of (A1)--(A2).} Since $f\in C^2_c(\R)$, $X^N_t$ is uniformly bounded for $t\geq 0$ and $N>1$, which trivially implies {\bf (A1)}.

To prove {\bf (A2)}, we will use the well-known duality between the voter model and coalescing random walks. More precisely, if we denote by $\{Y^{x,t}_{s}\}_{x \in \Z, t >0, s \leq t}$ a collection of coalescing random walks starting from each $x\in\Z$ at each time $t>0$, evolving backwards in time, each with increment distribution $p(\cdot)$, then there
exists a coupling between $(\eta_t)_{t\geq 0}$ and $\{Y^{x,t}_{s}\}_{x \in \Z, t >0, s \leq t}$ (using the so-called graphical construction) such that almost surely,
$$
\eta_{t}(x) = \eta_{0}(Y^{x,t}_{t}) \quad \mbox{for all } x \in \Z \mbox{ and } t >0.
$$
For more details, see e.g.~\cite{L85} or \cite{CD95}.

We start with a random walk estimate. Let $W:=(W_s)_{s\geq 0}$ be a continuous time random walk on $\Z$ with jump rate $2$
and jump kernel $p^*(\cdot)$, with $p^*(x)=\frac{p(x)+p(-x)}{2}$ for $x\in\Z$. Note that if $W_0=x-y$, then $(W_s)_{0\leq s\leq t}$ equals in law to
$(Y^{x,t}_{s} - Y^{y,t}_{s})_{0\leq s\leq t}$. Let
\be\label{tau}
\tau = \inf \{ s \geq 0 : W_{s} = 0 \}.
\ee
Let $\P_z(\cdot)$ and $\E_z[\cdot]$ denote respectively probability and expectation for $W$ with $W_0=z\in\Z$. Then for all $0 \leq s \leq t$ and $a >0$, we have
\be
\P_0( |W_{s}| \geq  a  ) \leq \frac{\E[W^{2}_{s}] }{a^{2}} = \frac{2\sigma^{2}s}{a^{2}}.
\ee
Using this bound and the strong Markov property, for all $z\in\Z$ and $a, s>0$, we have
\begin{eqnarray*}
\P_z(\tau\leq s, |W_s|\geq a) &=& \E_z[1_{\{\tau\leq s\}} \P_0(|W_{s-\tau}|\geq a)] \\
&\leq& \frac{2\sigma^2}{a^2} \E_z[1_{\{\tau\leq s\}} (s-\tau)] \\
&\leq& \frac{2\sigma^2 s}{a^2} \P_z(\tau\leq s).
\end{eqnarray*}
Therefore for $a>\sigma\sqrt{2s}$,
\be\label{hit1}
\P_z(\tau\leq s, |W_s|< a) \geq \Big(1-\frac{2\sigma^2 s}{a^2}\Big)\P_z(\tau\leq s).
\ee
On the other hand, for $a< |z|$, we have
\be\label{hit2}
\P_z(|W_s|<a) = \P_0(|W_s-z|<a) \leq \P_0\big(|W_s|\geq |z|-a\big) \leq \frac{2\sigma^2s}{(|z|-a)^2}.
\ee
Combining (\ref{hit1}) and (\ref{hit2}), and setting $a=4\sigma^2 s$ then gives
\be\label{coales}
\P_z(\tau\leq s) \leq \frac{2\sigma^2s}{(|z|-a)^2} \Big/ \Big(1-\frac{2\sigma^2 s}{a^2}\Big) = \frac{4\sigma^2s}{(|z|-2\sigma\sqrt{s})^2}
\ee
for all $s>0$ and $z\in\Z$ with $|z|>2\sigma \sqrt{s}$.
\bigskip

Let $T$, $\tau_N$ and $\delta_N$ be as in {\bf (A2)}. Let $M$ be such that the support of $f$ is contained in $[-M,M]$. Fix an $\delta >0$ small, and assume that $N$ is large enough so that $\delta_{N} < \delta$. Let $\E^{\zeta}[\cdot]$ denote expectation with respect to the voter model $(\eta_t)_{t\geq 0}$ with initial configuration
$\eta_0:=\zeta\in \{0,1\}^\Z$, and let ${\rm Var}^\zeta(\cdot)$ denote the corresponding variance. Then for any $\eps>0$,
\begin{eqnarray} \nonumber
\P( | X^{N}_{\tau_{N} + \delta_{N}} - X^{N}_{\tau_{N}} | > \eps) &\leq&
 \frac{1}{\eps^{2}} \E [ (X^{N}_{\tau_{N} + \delta_{N}} - X^{N}_{\tau_{N}})^{2}] =  \frac{1}{\eps^{2}} \E\big[\E^{\eta_{\tau_N N^2}}\big[ (X^{N}_{\delta_{N}} - X^{N}_0)^{2}\big]\big] \\
&\leq& \frac{2}{\eps^{2}}\E\Big[\big|\E^{\eta_{\tau_N N^2}}\big[X^{N}_{\delta_{N}} - X^{N}_0\big]\big|^2\Big] +
\frac{2}{\eps^{2}}\E\big[{\rm Var}^{\eta_{\tau_N N^2}}(X^{N}_{\delta_{N}})\big], \label{first}
\end{eqnarray}
where in the equality we used the strong Markov property for $(\eta_t)_{t\geq 0}$, and in the last inequality we added and subtracted
$\E^{\eta_{\tau_N N^2}}[X^{N}_{\delta_{N}}]$ from $X^{N}_{\delta_{N}} - X^{N}_0$ and used $(a+b)^2\leq 2(a^2+b^2)$.

Fix any $\zeta \in \{0,1\}^{\Z}$, and assume the coupling mentioned before between the voter model $(\eta_t)_{t\geq 0}$ with $\eta_0=\zeta$ and the collection of coalescing random walks $\{Y^{x,t}_{s}\}_{x \in \Z, t >0, s \leq t}$. Then
$$
\begin{aligned}
{\rm Var}^\zeta(X^N_{\delta_N}) &=  \frac{1}{N^2}\E^{\zeta} \Big[ \Big(\sum_{x \in \Z} f\Big(\frac{x}{N}\Big) \eta_{N^{2}\delta_{N}}(x)\Big)^2 \Big]  -\frac{1}{N^2} \E^{\zeta}\Big[\sum_{x \in \Z} f\Big(\frac{x}{N}\Big) \eta_{N^{2}\delta_{N}}(x)\Big]^{2} \\
&=  \frac{1}{N^{2}} \sum_{x,y \in \Z} f\Big(\frac{x}{N}\Big) f\Big(\frac{y}{N}\Big)
\Big(\E^{\zeta}\big[\eta_{N^{2}\delta_{N}}(x)  \eta_{N^{2}\delta_{N}}(y)\big]  -   \E^{\zeta}\big[\eta_{N^{2}\delta_{N}}(x)\big] \E^{\zeta}\big[\eta_{N^{2}\delta_{N}}(y)\big] \Big) \\
&=  \frac{1}{N^{2}} \sum_{x,y \in \Z} f\Big(\frac{x}{N}\Big) f\Big(\frac{y}{N}\Big) \Big(\E\big[\zeta(Y^{x,N^{2}\delta_{N}}_{N^{2}\delta_{N}})  \zeta (Y^{y,N^{2}\delta_{N}}_{N^{2}\delta_{N}})\big]  -   \E\big[\zeta (Y^{x,N^{2}\delta_{N}}_{N^{2}\delta_{N}})\big]  \E\big[  \zeta(Y^{y,N^{2}\delta_{N}}_{N^{2}\delta_{N}})\big] \Big) \\
&\leq \frac{|f|_\infty^2}{N^{2}}\!\!\!\!\!\!\sum_{x,y  \in [-MN,MN]\cap\Z} \!\!\!\!\!\!\!\!\!\!\!  \P(\tau^{x,y} \leq N^{2} \delta_{N})
\ \ \leq\ \ 4M|f|_\infty^2u+ \frac{|f|_\infty^2}{N^{2}} \!\!\!\!\!\!\sum_{x,y  \in [-MN,MN]\cap\Z\atop |x-y|\geq uN}\!\!\!\!\!\!\!\!\!\!\!  \P(\tau^{x,y} \leq N^{2} \delta_{N})
\end{aligned}
$$
for any $u>0$, where $\tau^{x,y} := \inf \{s \geq 0 :  Y^{x,N^{2}\delta_{N}}_{s}-Y^{y,N^{2}\delta_{N}}_{s}=0 \}$ is the time of coalescence of $Y^{x,N^{2}\delta_{N}}$ and $Y^{y,N^{2}\delta_{N}}$. Setting $u=4\sigma\delta^{\frac{1}{4}}$ and applying (\ref{coales}) with $s = N^{2}\delta_{N}$ and $z=x-y$ with $|x-y|\geq uN$
then gives
\be\label{bound1}
{\rm Var}^\zeta(X^N_{\delta_N}) \leq  C_1 \delta^{\frac{1}{4}} +  C_2  \delta^{\frac{1}{2}}
\ee
for some $C_1, C_2>0$ independent of $\zeta$ and all large $N$. This implies that the second term in (\ref{first}) is bounded by $2\eps^{-2}(C_1\delta^{\frac{1}{4}}+ C_2  \delta^{\frac{1}{2}})$. Since $\delta>0$ can be chosen arbitrarily small, this implies that the second term in (\ref{first}) tends to $0$ as $N\to\infty$.

To bound the first term in (\ref{first}), let $p_t(\cdot)$ denote the distribution of $Y^{x,t}_t-x$, which is the same for all $x\in\Z$. Then for any $\zeta\in \{0,1\}^\Z$,
\begin{eqnarray} \label{bound2}
 \big|  \E^{\zeta}[X^{N}_{\delta_{N}} - X^{N}_{0}] \big| &=& \Big| \E^{\zeta}\Big[ \frac{1}{N} \sum_{x \in \Z} f\Big(\frac{x}{N}\Big) \eta_{N^{2}\delta_{N}}(x) \Big] - \frac{1}{N} \sum_{x \in \Z} f\Big(\frac{x}{N}\Big) \zeta(x)\Big| \nonumber\\
 &=&  \frac{1}{N}\Big|\sum_{x \in \Z} f\Big(\frac{x}{N}\Big) \sum_{y \in \Z}  p_{N^{2}\delta_{N}}(y-x) \zeta(y)  - \sum_{y \in \Z} f\Big(\frac{y}{N}\Big)\zeta(y) \Big| \nonumber\\
  &=&   \frac{1}{N}\Big|  \sum_{x,y \in \Z}  p_{N^{2}\delta_{N}}(y-x)\Big(f\Big(\frac{x}{N}\Big) -f\Big(\frac{y}{N}\Big) \Big)\zeta(y)\Big| \nonumber\\
  &=&  \frac{1}{N} \Big|  \sum_{x,y \in \Z} p_{N^{2}\delta_{N}}(y-x)\Big(f^{\prime}\Big(\frac{y}{N}\Big)\frac{(x-y)}{N} + f^{''}(c_{N}(x,y))\frac{(x-y)^{2}}{2N^{2}} \Big)\zeta(y) \Big| \nonumber\\
    &\leq&  \frac{1}{N} \sum_{y \in \Z}  \sum_{x \in \Z}  |f^{\prime \prime}(c_{N}(x,y))| p_{N^{2}\delta_{N}}(y-x) \frac{(x-y)^{2}}{2N^{2}}, \nonumber
\end{eqnarray}
where we have applied Taylor expansion to $f(\frac{x}{N})$, for some $c_N(x,y)$ between $\frac{x}{N}$ and $\frac{y}{N}$, and in the inequality we have used the fact that $p_{N^{2}\delta_{N}}(\cdot)$ has zero mean. Since $|f''|_\infty<\infty$, $f''(c_N(x,y))\neq 0$ only if either $\frac{x}{N}$ or
$\frac{y}{N}$ is in the support of $f$,  and $p_{N^{2}\delta_{N}}(\cdot)$ has second moment $N^{2}\delta_{N}\sigma^2$, we see that the last term above
is bounded $C_3\delta_N$ for some $C_3>0$ independent of $\zeta$ and $N$. This implies that the first term in (\ref{first}) also tends to $0$ as $N\to\infty$.
The proof of {\bf (A2)} is then complete.
\qed

\subsection{Convergence of finite-dimensional distributions}
Let $\nu_t({\rm d}x):= 1_{\{x<\sigma B_t\}}{\rm d}x$ for a standard Brownian motion $(B_t)_{t\geq 0}$. By \cite[Prop.~11.1.VIII]{DV08}, the weak convergence of
$\mu^N_t$ to $\nu_t$ is equivalent to the weak convergence of $X^N_t(f)$ to $X_t(f):=\int f(x)\nu_t({\rm d}x)$ for every $f\in C_c(\R)$. Similarly, for any
$0\leq t_1< t_2\cdots < t_k$, the weak convergence of $(\mu^N_{t_1},\cdots, \mu^N_{t_k})$ to $(\nu_{t_1},\cdots, \nu_{t_k})$ is equivalent to the weak
convergence of
\be\label{fdd}
(X^N_{t_1}(f_1),\cdots, X^N_{t_k}(f_k)) \Asto{N} (X_{t_1}(f_1),\cdots, X_{t_k}(f_k)) \qquad \forall\, f_1,\cdots, f_k\in C_c(\R).
\ee
Since a.s.\ $|X_{t_i}(f_i)|\leq |f_i|_1$, (\ref{fdd}) would follow from the convergence of the moments, i.e.,
$$
\E\Big[\prod_{i=1}^k \big(X^N_{t_i}(f_i)\big)^{m_i}\Big] \asto{N} \E\Big[\prod_{i=1}^k \big(X_{t_i}(f_i)\big)^{m_i} \Big] \qquad \forall\, m_1,\cdots, m_k \in\N\cup\{0\}.
$$
Therefore (\ref{fdd}) will follow by showing that
\be\label{moment}
\E\Big[\prod_{i=1}^k X^N_{t_i}(f_i)\Big] \asto{N} \E\Big[\prod_{i=1}^k X_{t_i}(f_i)\Big] \quad \forall\, k\in\N, \, 0\leq t_0\leq \cdots \leq t_k, \, f_1,\cdots, f_k\in C_c(\R).
\ee
By the duality between $(\eta_t)_{t\geq 0}$ and the coalescing random walks $\{Y^{x,t}_{s}\}_{x \in \Z, t >0, s \leq t}$, we have
\begin{eqnarray*}
\E\Big[\prod_{i=1}^k X^N_{t_i}(f_i)\Big]
                  &=&  \frac{1}{N^k}\E \Big[\prod_{i=1}^{k} \Big(\sum_{x_{i} \in \Z} f_i\Big(\frac{x_{i}}{N}\Big) \eta_{N^{2}t_{i}}(x_{i})\Big) \Big]\\
                  &=&  \frac{1}{N^{k}}\sum_{u^N_{1},u^N_{2}, \ldots u^N_{k} \in \frac{\Z}{N}} \prod_{i=1}^k f_i(u^N_i)
                  \P\big(Y^{Nu^N_i,N^2t_i}_{N^2t_i}\leq 0 \mbox{ for all } 1\leq i\leq k\big).
\end{eqnarray*}
Firstly we note that as $N\to\infty$, the sequence of measures $\frac{1}{N^k} \sum_{u^N_{1},\ldots u^N_{k} \in \frac{\Z}{N}} \prod_{i=1}^k f_i(u^N_i) \delta_{u^N_i}(u_i)$ converges weakly to the finite measure $\prod_{i=1}^k f_i(u_i) {\rm d} u_i$ on $\R^k$. Secondly, it was shown in \cite[Section 5]{NRS05}
that if $(u^N_i, t_i)\to (u_i, t_i)$ for $1\leq i\leq k$, then
$$
\Big(\frac{1}{N}Y^{Nu^N_1,N^2t_1}_{sN^2}, \cdots, \frac{1}{N}Y^{Nu^N_k,N^2t_k}_{sN^2}\Big)_{s\geq 0} \Asto{N} (W^{u_1, t_1}_s, \cdots, W^{u_k, t_k}_s)_{s\geq 0},
$$
where $\Rightarrow$ denotes weak convergence, and $(W^{u_i,t_i}_\cdot)_{1\leq i\leq k}$ is a collection of backwards coalescing Brownian motions starting at $(u_i, t_i)_{1\leq i\leq k}$, each with diffusion coefficient $\sigma^2$. Only a finite second moment is required for such a convergence (actually only discrete time random walks were considered in~\cite{NRS05}, however the proof is easily seen to apply
to continuous time random walks as well). This implies that if $(u^N_i, t_i)_{1\leq i\leq k} \to (u_i, t_i)_{1\leq i\leq k}$, then
$$
\P\big(Y^{Nu^N_i,N^2t_i}_{N^2t_i}\leq 0 \mbox{ for all } 1\leq i\leq k\big) \asto{N} \P\big(W^{u_i, t_i}_{t_i} < 0 \mbox{ for all } 1\leq i\leq k\big).
$$
The above observations together imply that
\be\label{moment2}
\E\Big[\prod_{i=1}^k X^N_{t_i}(f_i)\Big] \asto{N} \idotsint \P\big(W^{u_i, t_i}_{t_i} < 0 \mbox{ for all } 1\leq i\leq k\big) \prod_{i=1}^k f_i(u_i) {\rm d}u_i.
\ee
We now recall that there is a natural coupling between forward and backward coalescing Brownian motions (see~\cite{STW00} and the later formulation in terms of the Brownian web and its dual in~\cite{FINR04, FINR06}). More specifically, there is a coupling between $(\sigma B_t)_{t\geq 0}$ for a standard Brownian motion $(B_t)_{t\geq 0}$ running forward in time, and $(W^{u_i,t_i}_\cdot)_{1\leq i\leq k}$ running backwards in time, such that $\sigma B_\cdot$ does not cross any $W^{u_i,t_i}_\cdot$. Therefore
$$
\P\big(W^{u_i, t_i}_{t_i} < 0 \mbox{ for all } 1\leq i\leq k\big) = \P(u_i < \sigma B_{t_i} \mbox{ for all } 1\leq i\leq k ).
$$
Substituting this identity into (\ref{moment2}) gives precisely (\ref{moment}).
\qed
\bigskip

\noindent
{\bf Acknowledgement} Both S.R.A.\ and R.S.\ are supported by grant R-146-000-119-133 from the National University of Singapore.

\end{document}